\documentclass{svproc}
\usepackage{url}

\usepackage{amssymb}
\usepackage{amsmath}
\usepackage{tikz}
\usepackage{mathtools}
\usepackage{stmaryrd}
\usepackage{todonotes}

\usetikzlibrary{patterns}
\usetikzlibrary{arrows.meta}

\newcommand{\lavg}{\left\lbrace\mskip-5mu\lbrace\mskip-2mu}
\newcommand{\ravg}{\right\rbrace\mskip-5mu\rbrace\mskip-2mu}

\begin{document}
\mainmatter              
\title{Domain of Dependence stabilization for the acoustic wave equation on 2D cut-cell meshes}
\titlerunning{DoD stabilization}  
%
\author{Gunnar Birke\inst{1} \and Christian Engwer\inst{1}\and Sandra May\inst{2} \and Florian Streitbürger\inst{3}}
\authorrunning{Gunnar Birke et al.} 
%
\tocauthor{Gunnar Birke, Christian Engwer, Sandra May}
\institute{Münster University, Applied Mathematics, 
\email{gunnar.birke@uni-muenster.de}
\and
Uppsala University, Department of Information Technology
\and
Dortmund University, Fakultät für Mathematik
}

\maketitle              

\begin{abstract}
Cut-cell meshes are an attractive alternative to avoid common mesh generation problems. For hyperbolic problems they pose additional challenges, as elements can become arbitrarily small, leading to prohibitive time step restrictions for explicit time stepping methods.
To alleviate this small cell problem we consider a particular stabilization method, the Domain of Dependence (DoD) method. So far, while posessing many favorable theoretical properties, in two dimensions the DoD method was essentially restricted to the transport equation. 
In this work we extend the DoD method to the acoustic wave equation in two dimensions and provide numerical results for validation.
\keywords{cut-cells, small cell problem, discontinuous Galerkin method, DoD stabilization, wave equation}
\end{abstract}

\section{Introduction}

The generation of traditional body fitted meshes can become very involved and time consuming, when working with complex geometries. One possible alternative is the use of cut-cell meshes. The drawback is that one has no control over the cut-cell shapes, in particular cut-cells can become arbitrarily small. Explicit time stepping methods, which are commonly used to simulate hyperbolic conservation laws, require a time step size that is chosen based on the smallest cell in the mesh. This becomes infeasible on a cut-cell mesh.
Instead, one likes to choose the time step based on the size of the larger uncut cells. This is referred to as the small cell problem.

As discontinuous Galerkin (dG) methods are attractive to solve hyperbolic partial differential equations (PDEs), different stabilization approaches have been developed in recent years to handle the small cell problem on cut-cell meshes, see for example \cite{giuliani}, \cite{kronbichler}. We consider the Domain of Dependence (DoD) stabilization method which was introduced in \cite{emns} for the linear transport equation and extended, in one space dimension, to non-linear equations in \cite{ms}. In two dimensions, the method was originally restricted to certain flow/geometry combinations, e. g. a flow parallel to a domain boundary. 
In \cite{PAMM} the stabilization was
generalized for $P^0$ discretizations of linear systems to handle triangular cut-cells with multiple inflow or outflow faces.


In this contribution we extend the lowest order DoD stabilization such that it can be applied to the acoustic wave equation. The main difference is that the wave equation cannot be written as a system of coupled transport equations, as the system is not globally diagonalizable. This adds significant complications in deriving the stabilization term. We focus on a $P^0$ dG discretization here, the development of appropriate higher-order stabilization terms for the acoustic wave equation is ongoing research.
The construction of the DoD stabilization went hand in hand with an $L^2$-stability analysis for the semi-discretization in space, see \cite{ms,ICOSAHOM,PAMM}, and the terms were designed such that we regain the spatial stability properties of the original dG scheme. This also holds true for the acoustic wave equation and the analysis was an important cornerstone in the development of the stabilization term. As the actual $L^2$-stability proof goes beyond the scope of this paper, we add appropriate remarks where necessary to understand the design of the actual formulation.

The outline of the paper is as follows: We first describe the problem setup and the underlying dG scheme. Then we introduce the DoD stabilization for the acoustic wave equation. We conclude with numerical results to validate our findings.

\section{Problem setup}

Let $\Omega \subset \mathbb{R}^2$ be an open domain.
We consider the acoustic wave equation for the solution $u = (p, v_1, v_2)$ given by
\begin{align*}
u_t + A_1 u_x + A_2 u_y & = 0 \:\:\: \text{ in } \Omega,\\
\tau u & = g \:\:\text{ on } \partial \Omega,\\
u(\cdot, 0) & = u_0 \:\:\text{ in } \Omega,
\end{align*}
with $p$ being the pressure, $v = (v_1, v_2)^t$ being the velocity and $\tau$ being a boundary operator imposing an inflow boundary condition for incoming waves. The system matrices are given as
\begin{displaymath}
A_1 = \footnotesize\begin{pmatrix}
0 & c & 0\\
c & 0 & 0\\
0 & 0 & 0
\end{pmatrix}\normalsize \text{ and }
A_2 = \footnotesize\begin{pmatrix}
0 & 0 & c\\
0 & 0 & 0\\
c & 0 & 0
\end{pmatrix}.
\end{displaymath}
Here $c > 0$ denotes the speed of sound. We want to stress that $A_1 A_2 \neq A_2 A_1$.

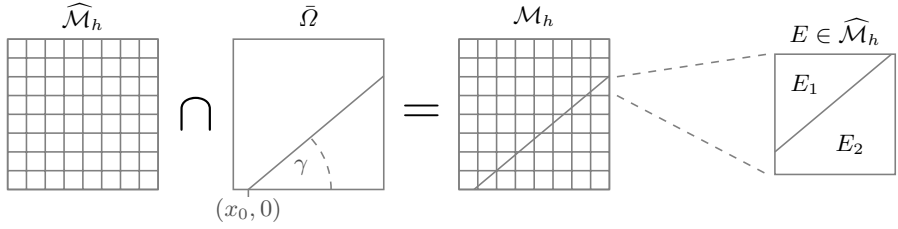
\begin{figure}
\centering
\begin{tikzpicture}[color=black,semithick,scale=1.0]
\begin{scope}
\node at (0.0,1.3) {\footnotesize${\widehat{\mathcal{M}}}_h$};
\draw[step=0.25cm,color=gray] (-1,-1) grid (1,1);
\draw[color=gray] (-1,-1) -- (1,-1) -- (1,1) -- (-1,1) -- cycle;
\node at (1.5,0.0) {\huge$\cap$};
\end{scope}
\begin{scope}[xshift=3.0cm]
\node at (0.0,1.3) {\footnotesize${\bar{\Omega}}$};
\draw[color=gray] (-1,-1) -- (1,-1) -- (1,1) -- (-1,1) -- cycle;
\draw[color=gray] 
(-0.7999, -1.1)
--
(-0.7999, -1.0) node[anchor=north]{\color{black!70}\footnotesize$(x_0, 0)$} -- (1.0, 0.5102);
\draw[color=gray, dashed] (0.3, -1.0) arc(0:44:1.0);
\node[color=black!70] at (-0.1,-0.7) {$\gamma$};
\node at (1.5,0.0) {\huge$=$};
\end{scope}
\begin{scope}[xshift=6.0cm]
\node at (0.0,1.3) {\footnotesize${\mathcal{M}_h}$};
\draw[step=0.25cm,color=gray] (-1,-1) grid (1,1);
\draw[color=gray] (-1,-1) -- (1,-1) -- (1,1) -- (-1,1) -- cycle;
\draw[color=gray] (-0.7999, -1.0) -- (1.0, 0.5102);
\end{scope}
\begin{scope}[xshift=10cm]
\draw[dashed,color=gray] (-2.9, 0.5) -- (-0.9,0.8);
\draw[dashed,color=gray] (-2.9, 0.25) -- (-0.9,-0.8);
\draw[color=gray] (-0.8,-0.8) -- (0.8,-0.8) -- (0.8,0.8) -- node[above]{\footnotesize\color{black}$E \in \widehat{\mathcal{M}}_h$} (-0.8,0.8) -- cycle;
\draw[color=gray] (-0.8,-0.5) -- (0.75,0.8);
\node at (-0.4,0.4) {\footnotesize $E_1$};
\node at (0.2,-0.4) {\footnotesize $E_2$};
\end{scope}
\end{tikzpicture}
\centering
\caption{Construction of the mesh: Out of the structured grid $\widehat{\mathcal{M}_h}$ on the domain $\Omega$ the mesh $\mathcal{M}_h$ is constructed by introducing cut-cells $E_1, E_2 \subset E \in \widehat{M_h}$ along the cut such that $\bar{E_1} \cup \bar{E_2} = \bar{E}$.}
\label{fig:cutted-square-geometry}
\end{figure}

In our numerical tests we choose $\Omega = [0, 1]^2$ and discretize it
by a structured grid $\widehat{\mathcal{M}_h}$. We then introduce an
artificial cut, a straight line going through the square, starting at
$(x_0, 0)$ and having an angle $\gamma$ relative to the $x$-axis. This
creates an internal boundary with two subdomains which we will resolve
by a cut-cell mesh $\mathcal{M}_h$. A sketch is contained in
Fig. \ref{fig:cutted-square-geometry}. The wave equation is then
solved in the whole domain $\Omega$.

Our discrete function space is defined as
\begin{displaymath}
\mathcal{V}_h = \mathcal{V}_h^0 = \{ v_h \in L^2(\Omega)^3 : (v_{h})_i\big|_E \in \mathcal{P}^0(E) \, \forall \, 1 \leq i \leq 3, \, E \in \mathcal{M}_h \}.
\end{displaymath}

We define the sets of internal and external faces as
\begin{align*}
\mathcal{F}^{\text{int}}_h & =  \{ F = \partial E_1 \cap \partial E_2 : E_1, E_2 \in \mathcal{M}_h, \; E_1 \neq E_2, \; |F| > 0 \},\\
\mathcal{F}^{\text{ext}}_h & =  \{ F = \partial E \cap \partial \Omega : E \in \mathcal{M}_h, \; |F| > 0 \}
\end{align*}
and the set of neighbor cells 
$
    \mathcal{N}(E) := \{ E' \in \mathcal{M}_h : |\bar{E}' \cap \bar{E}| > 0 \}.
$

So for any internal face $F \in \mathcal{F}^{\text{int}}_h$ there are always two unique elements $E_1, E_2 \in \mathcal{M}_h$ such that $\bar{E}_1 \cap \bar{E}_2 = F$. This face will often be denoted by $F_{E_1, E_2} = F_{E_2, E_1}$. We fix once and for all an orientation on $F = F_{E_1, E_2}$ by setting its outer normal vector $n$ to be $n = n_F = n_F(x) \coloneqq n_{E_1}(x)$ for $x \in F$ where $n_{E_1}$ is the outer unit normal field on $\partial E_1$. The flux matrix in normal direction on the face $F$ will be denoted by 
\[
A_F \coloneqq (n_F)_1 A_1 + (n_F)_2 A_2 = O_F \Lambda_F O_F^t
\]
where $O_F \Lambda_F O_F^t$ is an eigenvalue decomposition of the matrix $A_F$ with $\Lambda_F$ being a diagonal matrix and $O_F$ being an orthonormal matrix. Based on this, we define matrices which encode the flux directions as
\begin{align*}
  & A_F^+ = O_F \Lambda_A^+ O^t_F, \quad
    A_F^- = O_F \Lambda_F^- O^t_F,\quad\text{with }
    (\Lambda_F^{\pm})_{i,i} = \frac 1 2 \big(|(\Lambda_F)_{i,i}| \pm (\Lambda_F)_{i,i}\big).
\end{align*}
Note that $A_F = A_F^+ + A_F^-.$ We also introduce a generalization of the absolute value for such flux matrices by
$
|A_F| = A_F^+ - A_F^-.
$

An element $v_h \in \mathcal{V}_h$ is multi-valued on any internal face $F \in \mathcal{F}^{\text{int}}_h$. We define its average and jump by
\begin{displaymath}
\lavg v_h \ravg :=\frac{1}{2}(v_h\big|_{E_1}+v_h\big|_{E_2}), \quad
\left \llbracket v_h \right \rrbracket :=v_h\big|_{E_1}-v_h\big|_{E_2}.
\end{displaymath}
For exterior faces $F \in \mathcal{F}^{\text{ext}}_h$ we simply choose the unit outer normal and extend the definition of jump and average appropriately.

The (unstabilized) upwind semi-discretization in space is then given as: Find $u_h(t) \in \mathcal{V}_h$ such that
\begin{equation}
(\partial_t u_h(t), v_h)_{L^2(\Omega)} + a^{\text{upw}}_h(u_h(t), v_h) + l_h(v_h) = 0 \quad \forall \: v_h \in \mathcal{V}_h
\end{equation}
with
\begin{align*}
a^{\text{upw}}_h(u_h, v_h) & = \sum_{F \in \mathcal{F}_h^{\text{int}}} \int_F \langle A_F \left\lbrace\mskip-5mu\lbrace{u_h}\right\rbrace\mskip-5mu\rbrace, \left \llbracket{v_h}\right\rrbracket \rangle + \frac{1}{2} \langle |A_F| \left \llbracket{u_h}\right\rrbracket, \left \llbracket{v_h}\right\rrbracket \rangle ds\\
& + \sum_{F \in \mathcal{F}_h^{\text{ext}}} \int_F \langle A_F^+ u_h, v_h \rangle ds,\\
l_h(v_h) & = - \sum_{F \in \mathcal{F}_h^{\text{ext}}} \int_F \langle A_F^- g, v_h \rangle ds.
\end{align*}
Here, $\langle \cdot,\cdot \rangle$ denotes the standard scalar product in $l^2$.
{While we use an upwind flux here, we want to
  note that our following construction for the stabilization terms
  also holds for the Lax-Friedrichs flux, where some terms cancel. We skipped this for brevity.}
Discretization in time is then accomplished by the explicit Euler scheme. If the time step choice does not reflect the size of smaller cut-cells, this causes stability issues, which is why we need stabilization terms.

\section{Stabilization}

The main idea of the DoD stabilization is to extend the numerical domain of dependence of the neighbor cells of a small cut-cell by introducing additional numerical fluxes on the cut-cell boundary depending on the extrapolated solutions of the neighbor elements. 

Let $\mathcal{I} \subset \mathcal{M}_h$ be the set of cut-cells which are supposed to be stabilized (in practise this set will be chosen based on a volume fraction with respect to the background cells). For simplicity we assume that cells in $\mathcal{I}$ are not neighbors of each other. The $P^0$ stabilization fluxes are collected in a term
\begin{displaymath}J^0_h(u_h, v_h) = \sum_{E \in \mathcal{I}} J^{0, E}_h(u_h, v_h)
\end{displaymath}
and the stabilized scheme reads: Find $u_h(t) \in \mathcal{V}_h$ such that for any $v_h \in \mathcal{V}_h$
\begin{displaymath}
(\partial_t u_h(t), v_h)_{L^2(\Omega)} + a_h(u_h(t), v_h) + J^0_h(u_h(t), v_h) + l_h(v_h) = 0.
\end{displaymath}

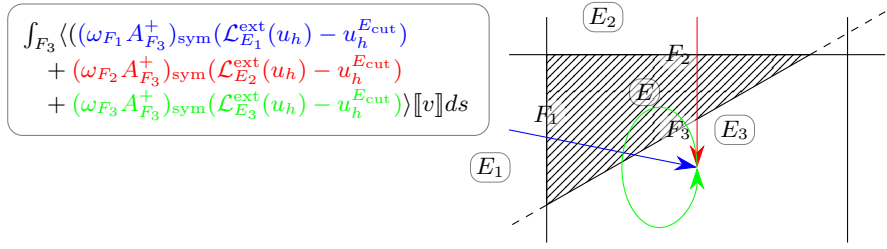
\begin{figure}[t]
\centering
\begin{tikzpicture}[whitebox/.style={%
        rounded corners=1.5mm,
        outer sep=0pt,
        inner sep=0.7mm,
        fill=white,
        opacity=0.5,
        text opacity=1},
        celltxt/.style={%
        rounded corners=1.5mm,
        outer sep=0pt,
        inner sep=0.7mm,
        draw, very thin,
        opacity=0.5,
        text opacity=1}
        ]
\begin{scope}[xshift=-1.0]
\fill[pattern=north east lines] (0.0, 2.0) -- (0.0, 0.0) -- (3.5, 2.0);
\draw[black] (0,0) 
  -- node[whitebox,shift={(0.0,0.2)}]{$F_1$}
     node[celltxt,shift={(-.75,-.5)}]{$E_1$} 
     ++(0,2)
  -- node[whitebox]{$F_2$}
     node[celltxt,shift={(-1,.5)}]{$E_2$} 
     ++(3.5,0)
  -- node[whitebox]{$F_3$}
     node[celltxt,shift={(.75,0)}]{$E_3$} 
     cycle;
\node[celltxt,whitebox] at (1.3, 1.5) {$E$};

\draw (4.5, 2) -- (-0.5, 2);
\draw (0, 2.5) -- (0, -0.5);
\draw (4, -0.5) -- (4, 2.5);
\draw[dashed] (3.5, 2.0) -- (4.5, 2.571);
\draw[dashed] (0, 0.0) -- (-0.5, -0.285);

\draw[-{Stealth[length=3mm, width=2mm]}, red](2.0,2.5) -- (2.0,0.5);
\draw[-{Stealth[length=3mm, width=2mm]}, blue](-0.5, 1.0) -- (2.0,0.5);
\draw[-{Stealth[length=3mm, width=2mm]},green] (2.0, 0.5) arc
    [
        start angle=0,
        end angle=360,
        x radius=0.5cm,
        y radius =0.8cm
    ] ;
\end{scope}
\node[
    rounded corners=2mm,
    outer sep=0pt,
    inner sep=2mm,
    draw=black!50
] at (-4, 1.8) {\parbox{6cm}{$\int_{F_3} \langle( \textcolor{blue}{(\omega_{F_1}A_{F_3}^+)_{\text{sym}}(\mathcal{L}^{\text{ext}}_{E_1}(u_h) - u_h^{E_{\text{cut}}})}\\
\phantom{+}+\textcolor{red}{(\omega_{F_2}A_{F_3}^+)_{\text{sym}} (\mathcal{L}^{\text{ext}}_{E_2}(u_h) - u_h^{E_{\text{cut}}})}\\
\phantom{+}+\textcolor{green}{(\omega_{F_3}A_{F_3}^+)_{\text{sym}} (\mathcal{L}^{\text{ext}}_{E_3}(u_h) - u_h^{E_{\text{cut}}})}\rangle \llbracket v \rrbracket ds$}};
\end{tikzpicture}
\caption{DoD stabilization for the acoustic equation on a triangular cut-cell. On the left are the stabilization fluxes on the face $F_3$ (without the parameter $\eta_E$), colored in correspondence with the arrows on the right indicating the different couplings. Note that there will be additional stabilization fluxes on $F_1$ and $F_2$. The green flux is somewhat of a curiosity and might just be a product of the way we construct our weighting matrices for the specific case of the acoustic equation. In a sense it acts similar to a reflecting wall.}
\label{fig:small-cut-cell-stabilization}
\end{figure}

For any $E \in \mathcal{M}_h$ we introduce a map $\mathcal{L}^{\text{ext}}_E: \mathcal{V}_h(E) \to \mathcal{P}^0(\Omega)^3$, called an extension operator, such that $\mathcal{L}^{\text{ext}}_E(u_h)\big|_E = u_h\big|_E$. The cell stabilization terms $J^{0, E}_h(u_h,v_h)$ then contain extended fluxes
\begin{displaymath}
\eta_E\!\! \sum_{(E_1, E_2) \in \mathcal{N}(E)} \int_{F_{E, E_2}} \!\!\Big\langle (\omega_{F_{E, E_1}} A_{F_{E, E_2}}^+)_{\text{sym}} (\mathcal{L}^{\text{ext}}_{E_1}(u_h) - u_h\big|_E), \llbracket v_h \rrbracket \Big\rangle ds,
\end{displaymath}
with appropriate weighting matrices $\omega_{F_{E, E_1}} \in \mathbb{R}^{3\times3}$ and a penalty parameter $\eta_E \in (0, 1]$. Here, we have denoted the symmetrization $\frac{1}{2}(A + A^t)$ for a matrix $A \in \mathbb{R}^{3 \times 3}$ by $A_{\text{sym}}$.
This introduces a direct mass transport from $E_1$ to $E_2$ for any pair $(E_1, E_2)$ of neighbors of a small cut-cell $E$. (For a triangular cut-cell, the sum contains 9 terms.) This is a result of the inherent nature of the wave equation, creating waves in all possible directions. The weighting matrices describe how much of an inflow coming from $E_1$ is transported to $E_2$. Figure \ref{fig:small-cut-cell-stabilization} shows an illustration. Similarly to what was proposed in \cite{PAMM}, we require for the weighting matrices that
\begin{subequations}
\begin{align}
\label{eq:flux-distribution} \sum_{E_2 \in \mathcal{N}(E)} \int_{F_{E, E_2}} \hspace*{-1em}\omega_{E, E_1} A^+_{F_{E, E_2}} ds & = - \int_{F_{E, E_1}} \hspace*{-1em} A^-_{F_{E, E_1}} ds \qquad \forall \: E_1 \in \mathcal{N}(E),\\
\label{eq:flux-normalization} \sum_{E_1 \in \mathcal{N}(E)} \omega_{F_{E, E_1}} & = \, \text{Id}_{3 \times 3}.
\end{align}
\end{subequations}
The first equation describes how incoming flow is distributed among the cut-cell's neighbors while the second equation ensures that the overall amount of flow over a single face is preserved.

In \cite{PAMM} we were able to prove a discrete dissipation and $L^2$-stability result. One key step of the proof is an application of the binomial formula relying on the symmetry of the matrices $\omega_{F_{E, E_1}} A_{F_{E, E_2}}^+$, which holds given simultaneously diagonalizable (equivalently, commuting) system matrices. 
For the acoustic equation this matrix product will in general not be symmetric and the proof breaks down. A working fix is to instead take the symmetrization of the matrix $\omega_{F_{E, E_1}} A_{F_{E, E_2}}^+$.


Taking the symmetrization is not enough by itself though since the resulting flux matrices potentially possess negative eigenvalues leading to incorrect fluxes. This can be corrected by introducing terms
\begin{displaymath}
-\kappa \eta_E\!\! \sum_{(E_1, E_2) \in \mathcal{N}(E)} \int_{F_{E, E_2}} \!\!\Big\langle (\omega_{F_{E, E_1}} A_{F_{E, E_2}}^+)_{\text{sym}}^- (\mathcal{L}^{\text{ext}}_{E_1}(u_h) - u_h\big|_{E_2}), \mathcal{L}^{\text{ext}}_{E_1}(v_h) - v_h\big|_{E_2} \Big\rangle ds.
\end{displaymath}
For any $\kappa \geq 1$ this term ensures
$L^2$-stability of the semi-discrete form. We observe that certain choices of $\kappa$ significantly improve the numerical approximation; the 
optimal choice still needs to be
investigated.
Note that this new term has an opposite sign, uses the negative part
of the symmetrized form, and subtracts the solution from cell
neighbor $E_2$ (instead of the cut-cell $E$).


Putting those components together leads to a cell stabilization term
\begin{align*}
J^{0, E}_h &(u_h, v_h) \\ & = \eta_E\!\! \sum_{(E_1, E_2) \in \mathcal{N}(E)} \int_{F_{E, E_2}} \!\!\Big\langle (\omega_{F_{E, E_1}} A_{F_{E, E_2}}^+)_{\text{sym}} (\mathcal{L}^{\text{ext}}_{E_1}(u_h) - u_h\big|_E), \llbracket v_h \rrbracket \Big\rangle ds\\
& -\kappa \eta_E\!\! \sum_{(E_1, E_2) \in \mathcal{N}(E)} \int_{F_{E, E_2}} \!\!\Big\langle (\omega_{F_{E, E_1}} A_{F_{E, E_2}}^+)_{\text{sym}}^- (\mathcal{L}^{\text{ext}}_{E_1}(u_h) - u_h\big|_{E_2}), \mathcal{L}^{\text{ext}}_{E_1}(v_h) - v_h\big|_{E_2} \Big\rangle ds.
\end{align*}
For this formulation together with the properties \eqref{eq:flux-distribution} and \eqref{eq:flux-normalization} one can prove $L^2$-stability for the semi-discretization in space. Unfortunately the proof is too long for this contribution to be included.

The weighting matrices can be directly computed from \eqref{eq:flux-distribution}. Due to the system matrices being divergence-free, this equation is equivalent to
\begin{displaymath}
\omega_{F_{E, E_1}} \Big(\sum_{E_2 \in \mathcal{N}(E)} |F_{E, E_2}| A_{F_{E, E_2}}^-\Big) = |F_{E, E_1}| A_{F_{E, E_1}}^-.
\end{displaymath}
The matrix $(\sum_{E_2 \in \mathcal{N}(E)} |F_{E, E_2}| A_{F_{E,E_2}}^-)$ is invertible since on a triangular cut-cell we got three summands, each of rank one, and the three column spaces are linearly independent. Thus we have
\begin{displaymath}
\omega_{F_{E, E_1}} = |F_{E, E_1}| A_{F_{E, E_1}}^- \Big(\sum_{E_2 \in \mathcal{N}(E)} |F_{E, E_2}| A_{F_{E, E_2}}^-\Big)^{-1}.  
\end{displaymath}

\section{Numerical results}

We present numerical results to support our findings. Our implementation uses the DUNE framework \cite{dune3,dune2}, in particular the dune-udg module \cite{dune-udg,dune-udg2} and the TPMC library \cite{tpmc}.
We consider an analytic test case \cite{example} with
\begin{displaymath}
u(x, y, t) = \frac 1 c\begin{pmatrix}
  - 
  \cos(2 \pi ct)(\sin(2\pi x) + \sin(2 \pi y)) \\
  \sin(2 \pi ct) \cos(2 \pi x) \\
  \sin(2 \pi ct) \cos(2 \pi y)
\end{pmatrix}
\end{displaymath}
Initial condition $u_0$ and inflow boundary conditions $g$ are given
by the exact solution. The speed of sound is chosen as $c = \frac{1}{2}$, the final time is $T = 0.3$. For the CFL-condition we choose $\Delta t = 0.3 \frac{\Delta x}{c}$ where $\Delta x = \frac{1}{N}$ with $N$ being the number of
fundamental cells in one direction of the grid. We choose values for
$N$ from the set $\{400, .., 1200\}$, starting at $400$ and with a
stepsize of $23$. In our experiments the volume fractions of cut-cells are in the interval $[7.24 \cdot 10^{-10}, 5.37 \cdot 10^{-5}]$. For the cut angle we set $\gamma = 35^{\circ}$ and choose $x_0 = 0.2001$ for the start of the cut. We classify a cell $E \in \mathcal{M}_h$ as a small cut-cell if $\frac{|E|}{(1/N)^2} \leq 0.4$. The penalty parameter $\eta_E$ is chosen as $(1 - \eta_E) = \frac{|E|}{ \Delta t c \max_{F \subset \partial E} |F|}$.

We compute the $L^2$-error of the discrete solution at the final time $T$ with respect to the exact solution. We also compute the pointwise error at certain quadrature points, which gives an approximation to the error in the $L^{\infty}$-norm. We will display results for $\kappa \in \{1.0, 7.5\}$ to display the effect of this additional parameter.

\begin{figure}[th]
\includegraphics[width=0.5\linewidth]{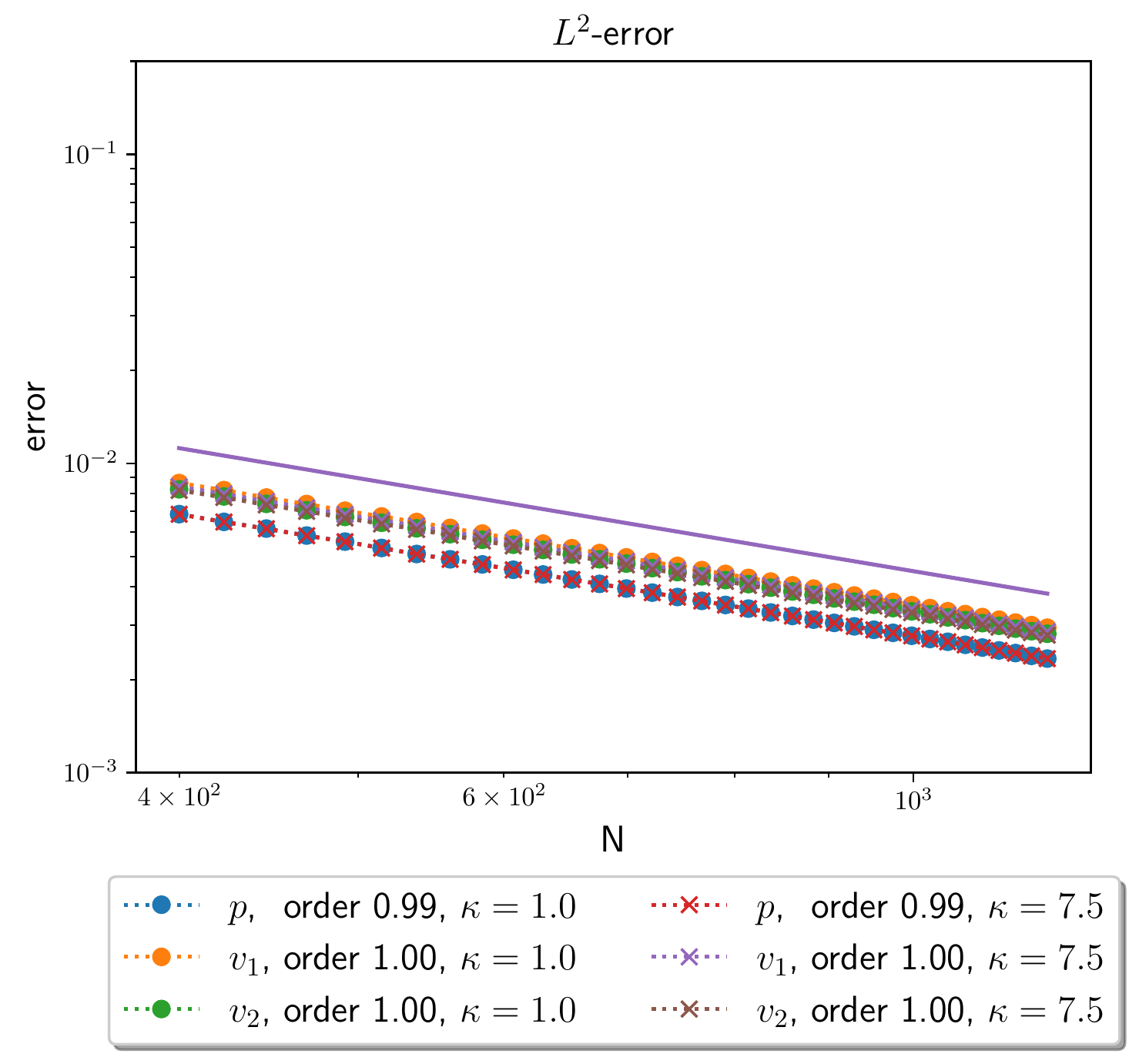}\hspace{-0.5em}
\includegraphics[width=0.5\linewidth]{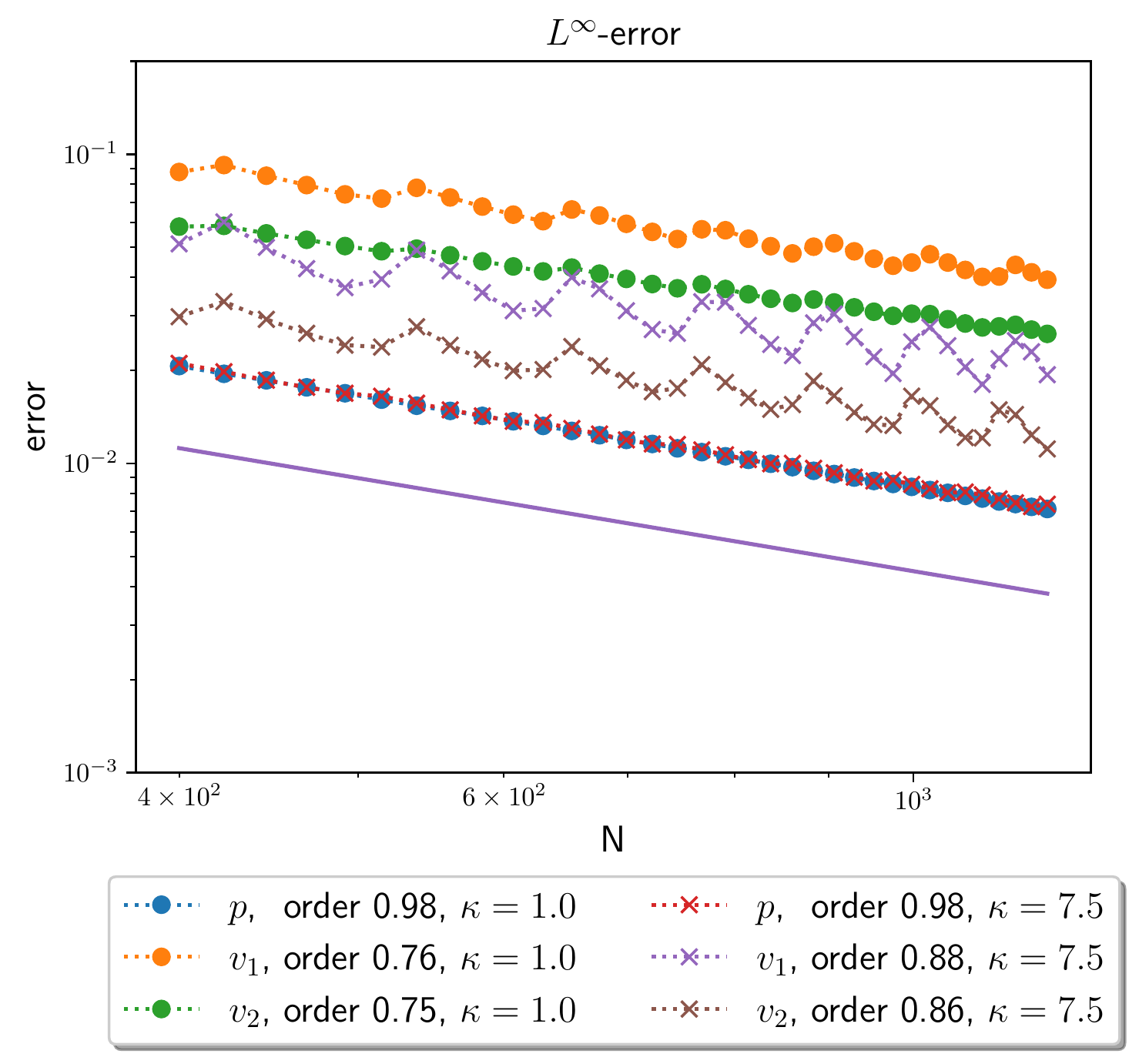}
\caption{Error at $T = 0.3$ in the $L^2$ norm (left) and in the $L^{\infty}$ norm (right). The straight purple line is for reference and denotes convergence of order 1.}
\label{fig:errors}
\end{figure}

The computed errors are plotted in figure \ref{fig:errors}. We observe the expected convergence order in the $L^2$-norm in all components, regardless of the choice of $\kappa$. For the $L^{\infty}$-norm the situation differs. The pressure component shows again optimal convergence behavior, independent of $\kappa$. However for the velocity components the convergence order is influenced by the choice of $\kappa$ and we obtain better absolute errors and a better convergence order (at the cost of some mild wiggles) for a larger $\kappa = 7.5$. The wiggles suggest that the parameter $\kappa$ should be chosen cell dependent, but how to choose it exactly requires further investigation.

\section{Discussion and Outlook}

We presented an extension of the DoD stabilization to the acoustic wave equation for $P^0$ trial and test functions on a structured grid with triangular cut-cells. We observed a decent convergence behavior in our numerical test setup. We note that the new parameter $\kappa$ included in the extension has a considerable influence on the measured error. Its precise effect and optimal choice will be investigated more in the future. An extension to higher-order approximations, as well as to non-linear systems, e. g. the Euler equations, is ongoing research.

\subsubsection{Acknowledgements} 
The authors acknowledge support by the Deutsche
For\-schungsgemeinschaft as project 439956613 under contract numbers
EN\,1042/5-1 and MA\,7773/4-1/2, as well as under Germany’s Excellence
Strategy EXC 2044\\ 390685587, Mathematics M\"unster: Dynamics --
Geometry -- Structure.
\end{document}